\documentclass[12pt]{amsart}
\usepackage{amssymb}
\usepackage{amsmath}
\usepackage{amsthm}
\usepackage{times}
\usepackage{graphicx}

\newtheorem{theorem}{Theorem}[section]
\newtheorem{lemma}[theorem]{Lemma}
\newtheorem{corollary}[theorem]{Corollary}
\newtheorem{proposition}[theorem]{Proposition}
\newtheorem{prop}[theorem]{Proposition}

\theoremstyle{definition}
\newtheorem{definition}[theorem]{Definition}

\newtheorem{question}[theorem]{Question}

\theoremstyle{remark}
\newtheorem{remark}[theorem]{Remark}

\numberwithin{equation}{section}

\newcommand{\calO}{\mathcal{O}}
\newcommand{\calL}{\mathcal{L}}

\newcommand{\ignore}[1]{}

\def\deg{{\text{deg}}}

\def\spmapright#1{\smash{%
   \mathop{\hbox to 1.3cm{\rightarrowfill}}
       \limits^{#1}}}

\def\sbmapright#1{\smash{%
   \mathop{\hbox to 1.3cm{\rightarrowfill}}
       \limits_{#1}}}

\newcommand{\mapright}[1]{%
\smash{\mathop{%
   \hbox to 1cm{\rightarrowfill}}\limits^{#1}}}

\newcommand{\mapleft}[1]{%
\smash{\mathop{%
   \hbox to 1cm{\leftarrowfill}}\limits^{#1}}}

\begin{document}

\title{On the multicanonical systems of \\
quasi-elliptic surfaces}
\author{Toshiyuki Katsura}
\thanks{Research of the first author is partially supported by JSPS Grant-in-Aid 
for Scientific Research (C) No. 20K03530, and the second author by (C) No. 17K05208} 

\address{Graduate School of Mathematical Sciences, The University of Tokyo, 
Meguro-ku, Tokyo 153-8914, Japan}
\email{tkatsura@ms.u-tokyo.ac.jp}

\author{Natsuo Saito}
\address{Graduate School of Information Sciences, Hiroshima City University,
Asaminami-ku, Hiroshima, 731-3194, Japan}
\email{natsuo@math.info.hiroshima-cu.ac.jp}


\begin{abstract}
We consider the multicanonical systems $\vert mK_{S}\vert$ of quasi-elliptic 
surfaces with Kodaira dimension $1$ in characteristic 2. We show that 
for any $m \geq 6$ 
$\vert mK_{S}\vert$ 
gives the structure 
of quasi-elliptic fiber space, and 6 is the best possible number to give the structure for any such surfaces.
\end{abstract}



\maketitle

\section{Introduction}
Let $k$ be an algebraically closed field of characteristic $p \geq 0$, and
let $S$ be a nonsingular complete algebraic surface 
with Kodaira dimension 1 defined over $k$.
Then, $S$ has a structure of genus 1 fibration $\varphi : S \longrightarrow B$.
We denote by $K_{S}$ a canonical divisor of $S$ and 
we consider the multicanonical system $\vert mK_{S}\vert$.
As is well known, the multicanonical system $\vert mK_{S}\vert$
gives the genus 1 fibration if $m$ is large enough.
In Katsura and Ueno \cite{KU} and Katsura \cite{K} (see also Iitaka \cite{I}), 
we considered 
the following question:

\begin{question}
(1) Does there exist a positive integer $M$ such that if $m\geq M$,
the multicanonical system $\vert mK_{S}\vert$ gives a structure
of genus 1 fibration for any elliptic 
surface $S$ over $k$ with Kodaira dimension $1$?

(2) What is the smallest $M$ which satisfies this property?
\end{question}

For this question, we have the following theorem.

\begin{theorem}
$(1)$ For the complex analytic elliptic surfaces, $M = 86$
and 86 is best possible (cf. Iitaka \cite{I}). 

$(2)$ For the algebraic elliptic surfaces, if the characteristic 
$p = 0$ or $p \geq 3$, 
then $M = 14$ and 14 is best possible 
(Katsura and Ueno \cite{KU} and Katsura \cite{K}). 

$(3)$ For the algebraic elliptic surfaces, if the characteristic $p =2$, then  
$M = 12$ and 12 is best possible (Katsura \cite{K}).
\end{theorem}

If $p = 2$ or $3$, there are two kinds of
genus 1 fibrations, namely, the elliptic fibration 
and the quasi-elliptic fibration (cf. Bombieri and Mumford \cite{BM}).
In these cases, we can also consider the same question for quasi-elliptic surfaces 
with Kodaira dimension 1.
In characteristic 3, we already showed the following results 
(Katsura \cite{K2}).

\begin{theorem}
For the quasi-elliptic surfaces in characteristic 3, we have $M = 5$, and $5$ is
best possible.
\end{theorem}

Therefore, the remaining case of the question 
for the surfaces with Kodaira dimension 1
is the one in characteristic 2, and
in this paper
we show the following theorem. It finishes the answer to the question above
for surfaces with Kodaira dimension 1
which S. Iitaka considered in the case of complex analytic elliptic surfaces
in 1970 (cf. \cite{I}).
\begin{theorem}
For the quasi-elliptic surfaces in characteristic 2,
we have $M = 6$ and 6 is best possible. 
\end{theorem}

In Section 2, we summarize basic facts on the theory of vector fields
in positive characteristic and some results on quasi-elliptic surfaces.
In Section 3, we give a criterion for a vector field 
that makes a singularity on the quotient of curve.
In Section 4, we construct a quasi-elliptic surface
over an elliptic curve with only one tame multiple fiber and 
examine the structure of its multicanonical system.
In Section 5, we examine the multicanonical system of
quasi-elliptic surfaces in characteristic 2 and show our main theorem.

The first author thanks Professor K. Ueno for indicating him
the original question and for many discussions.

\section{Preliminaries}
Let $k$ be an algebraically closed field of characteristic $p > 0$,
and let $S$ be a nonsingular complete algebraic surface defined over $k$.
A non-zero rational vector field $D$ on $S$ is called $p$-closed 
if there exists
a rational function $f$ on $S$ such that $D^p = fD$. 

We use a vector field to construct a quotient surface of $S$.
Let $\{U_{i} = {\rm Spec} A_{i}\}$ be an affine open covering of $S$ and we set 
$A_{i}^{D} = \{D(\alpha) = 0 \vert \alpha \in A_{i}\}$. 
Then, affine surfaces $\{U_{i}^{D} = {\rm Spec} A_{i}^{D}\}$ glue together to 
define a normal quotient surface $S^{D}$.

We now recall some results on vector fields by 
Rudakov and Shafarevich \cite[Section 1]{RS}.
Now, we  assume that $D$ is $p$-closed. Then, we know that
the natural morphism $\pi : S \longrightarrow S^D$ is a purely
inseparable morphism of degree $p$. 
If the affine open covering $\{U_{i}\}$ of $S$ is fine enough, then
taking local coordinates $x_{i}, y_{i}$
on $U_{i}$, we see that there exist $f_{i}, g_{i}\in A_{i}$ and 
a rational function $h_{i}$
such that the divisors defined by $f_{i} = 0$ and by $g_{i} = 0$ have 
no common divisor
and that the vector field $D$ is expressed as
$$
 D = h_{i}\left(f_{i}\frac{\partial}{\partial x_{i}} + g_{i}\frac{\partial}{\partial y_{i}}\right)
\quad \mbox{on}~U_{i}.
$$
Divisors $(h_{i})$ on $U_{i}$
give a global divisor $(D)$ on $S$, and zero-cycles defined
by the ideal $(f_{i}, g_{i})$ on $U_{i}$ give a global zero cycle 
$\langle D \rangle $ on $S$. A point contained in the support of
$\langle D \rangle $ is called an isolated singular point of $D$.
(\cite[Theorem 1, Corollary]{RS}). 
Rudakov and Shafarevich showed that $S^D$ is nonsingular
if and only if $\langle D \rangle  = 0$.
When $S^D$ is nonsingular,
they also showed a canonical divisor formula
\begin{equation}\label{canonical}
K_{S} \sim \pi^{*}K_{S^D} + (p - 1)(D),
\end{equation}
where $\sim$ means linear equivalence.

Now, we consider an irreducible curve $C$ on $S$ and we set $C' = \pi (C)$.
Take an affine open set $U_{i}$ above such that $C \cap U_{i}$ is non-empty.
The curve $C$ is said to be integral with respect to the vector field $D$
if $D$ is tangent to $C$ at a general point of $C \cap U_{i}$. 
Rudakov-Shafarevich
showed the following proposition (cf. \cite[Proposition 1]{RS}):

\begin{prop}\label{insep}

$({\rm i})$  If $C$ is integral, then $C = \pi^{-1}(C')$ and $C^2 = pC'^2$.

$({\rm ii})$  If $C$ is not integral, then $pC = \pi^{-1}(C')$ and $pC^2 = C'^2$.
\end{prop}

Now, let $\varphi : S \longrightarrow B$ be a quasi-elliptic surface.
We denote by $g$ the genus of the curve $B$.
As was shown in Katsura \cite{K2}, we have $Alb(S) \cong J(B)$, and 
$\chi ({\calO}_{S}) \geq (1 -g)/3$
(See also Lang~\cite{L} and Raynaud~\cite{R}). Here, $Alb(S)$
is the Albanese variety of $S$ and $J(B)$ is the Jacobian variety of $B$.
As a corollary, we know that
if $g = 1$, then $\chi ({\calO}_{S}) \geq 0$, and that
if $g = 0$, then $\chi ({\calO}_{S}) \geq 1$. We will freely use
these inequalities in Section 5.

\section{Cuspidal points}
From here on, let $k$ be an algebraically closed field of characteristic 2,
if otherwise mentioned.
Let $S$ be a nonsingular complete algebraic surface over $k$,
and let $D$ be a non-zero 2-closed rational vector field on $S$.
Let $U$ be an affine open set of $S$, and $x, y$ be local coordinates
of $U$. Then, as in Section 2, $D$ is given by 
$$
     D = h(f\partial /\partial{x} + g\partial /\partial{y}),
$$
where $f, g$ are regular functions on $U$ such that $f = 0$ has no
common curves with $g = 0$, and where $h$ is a rational function on $S$.

\begin{lemma}
Under the assumption above, $D(fg) = 0$ holds.
\end{lemma}
\proof{We set $\alpha = hf$ and $\beta = hg$. Since there exists
a rational function $\gamma$ such that $D^{2} = \gamma D$,
we have
$$
\begin{array}{l}
 \alpha \alpha_{x} + \beta \alpha_{y} = \gamma \alpha, \\
\alpha \beta_{x} + \beta \beta_{y} = \gamma \beta.
\end{array}
$$
Therefore, by direct calculation, we have $D(\alpha \beta) = 0$.
Since $\alpha \beta = h^{2}fg$, we conclude $D(fg) = 0$.
\hfill $\Box$}

\begin{corollary}\label{cor:0}
$D(f/g) = 0$.
\end{corollary}
\proof{We have $D(f/g) = D(fg/g^{2}) = (1/g^{2})D(fg) = 0$.
\hfill $\Box$}

\begin{definition}\label{cusp}
Let $D$ be a non-zero rational vector field on a nonsingular surface $S$,
and $C$ be a nonsingular irreducible curve on $S$. Let $P$ be a point
on $C$ which is not an isolated singular point of $D$. 
If $D$ is non-integral
on $C$ and integral at a point $P$ on $C$, we call $P$ a cuspidal
point of the vector field $D$.
\end{definition}

\begin{proposition}\label{cuspidal curve}
Under the notation in Definition \ref{cusp}, we consider 
the projection $\pi : S \longrightarrow S^{D}$.
Then, the image $\pi (P)$ of the cuspidal point $P$
is a singular point of the curve $\pi (C)$.
\end{proposition}
\proof{
Let ${O}_P$ be the local ring of the cuspidal point $P$
and let $x, y$ be a system of parameters of $O_P$.
Let $x = 0$ be a local equation of $C$ at the point $P$.
By the definition of cuspidal points, there exist elements
$\alpha$, $\beta$, $\gamma$ and $\delta$ of ${O}_P$ 
and a constant 
$c \in k$ such that $\beta \not\equiv 0$ and $c \neq 0$, and
such that $f = \alpha x + \beta y$ and $g = \gamma x + \delta y + c$.
Since the situation is local, we may omit $h$ from $D$.
By Corollary \ref{cor:0}, we see that $D(x + (f/g)y) = 0$.
Since $g(P) \neq 0$, $x + (f/g)y$ is contained in ${O}_P$.
Considering the completion ${\hat{O}}_P$ of $O_P$, we have 
${\hat{O}}_P \cong k[[x, y]]$.
Since $k[[x, y]]^{D} \supset k[[y^2, x + (f/g)y]]$ and
$\dim_k k[[x, y]]^{D}/k[[x^2, y^2]] = \dim_k k[[y^2, x + (f/g)y]]/k[[x^2, y^2]] = 2$,
we have $k[[x, y]]^{D} = k[[y^2, x + (f/g)y]]$. Although by the general theory
of the vector field the point $\pi (P)$ is a nonsingular point of $S^D$,
this result also shows that $S^D$ is nonsingular at $\pi (P)$.
We set $X = x^2$, $Y = y^2$ and
$Z = x + (f/g)y$, and let $\tilde{f}$, $\tilde{g}$ be elements of ${O}_P$
whose coefficients are the squares of the ones of $f$, $g$, respectively.
Let $S'$ be a surface defined by the equation
$$
Z^2 = X + (\tilde{f}/\tilde{g})Y.
$$
Since the degrees of the algebraic extensions $k(S)/k(S^D)$ and $k(S)/k(S')$ 
of fields are 2 and $k(S^D) \supset k(S')$ holds, we have $k(S^D) = k(S')$, 
that is,
$S'$ is birationally equivalent to $S^D$. Since $\tilde{g}(P)= c^2 \neq 0$, 
$S'$ is nonsingular at the point $(X, Y, Z) = (0, 0, 0)$. 
Therefore, by the Zariski main theorem
the surface $S^D$ is isomorphic to $S'$ around $\pi (P)$.  
The curve $\pi (C)$ is defined by
$X =0$ at the point $\pi (P)$.
Therefore, the equation of the curve $\pi (C)$ at $\pi (P)$ on the plane
$X = 0$ is given by
$$
     Z^2(\tilde{\delta}\vert_{X = 0}Y + c^2) = \tilde{\beta}\vert_{X = 0}Y^2.
$$
Here, the notation of $\tilde{\beta}$ and $\tilde{\delta}$ are 
similar to $\tilde{f}$ and $\tilde{g}$.
This equation for the curve $\pi (C)$ shows 
that $\pi (P)$ is a singular point of $\pi (C)$.
\hfill $\Box$}

\section{A construction of a quasi-elliptic surface}
Let $E$ be an elliptic curve and $\{U_0, U_{\infty}\}$
be an affine open covering and let $U_0$ (resp. $U_{\infty}$)
be given by the equation
$$
    y^2 + y = x^3 ~({\rm resp.}~z^2 + z = w^3).
$$
The change of coordinates is given by
$$
y = 1/z,~x = w/z.
$$
Let $\{V_0, V_{\infty}\}$ ($V_0 \cong V_{\infty}\cong {\bf A}^1$ : an affine line)
be affine open covering of the projective
line ${\bf P}^{1}$ and $t$ (resp.$s$) be a coordinate of $V_0$
(resp. $V_{\infty}$). The change of coordinates is given by
$$
   t = 1/s.
$$
We consider the algebraic surface $S = E \times {\bf P}^1$.
Then, $\{U_i\times V_j~\vert ~ i = 0, \infty; j = 0, \infty\}$
gives an affine open covering of $S$.  
We have a projection
$$
\psi : S \longrightarrow E.
$$
Let $C_{\infty}$ be the curve on $S$ defined by $s = 0$.
We consider the following rational vector field $D$ on $U_0\times V_0$.
$$
(I)\quad  D = y\frac{\partial}{\partial x} + 
 (x^2 + x^2t + t^4)\frac{\partial}{\partial t}.
$$
Then, $D$ gives a rational vector field on $S$ and on each affine chart
it is concretely given as follows:
$$
\begin{array}{lll}
(II)&D &= \frac{1}{z^2}\{z\frac{\partial}{\partial w} + 
 (w^2 + w^2t + z^2t^4)\frac{\partial}{\partial t}\} \\
    &  &= \frac{1}{w^4}\{(z+1)w\frac{\partial}{\partial w} + 
 ((z + 1)^2 + (z + 1)^2t + w^4t^4)\frac{\partial}{\partial t}\} \\
     & &\quad \quad \mbox{on}~U_{\infty}\times V_0 \\
(III)& D &= \frac{1}{s^2}\{ys^2\frac{\partial}{\partial x} + 
 (x^2s^4 + x^2s^3 + 1)\frac{\partial}{\partial s}\}\\
  & &\quad \quad   \mbox{on}~U_0\times V_{\infty}\\
(IV) & D &= \frac{1}{z^2s^2}\{zs^2\frac{\partial}{\partial w} + 
 (w^2s^4 + w^2s^3 + z^2)\frac{\partial}{\partial s}\} \\
 &  &  = \frac{1}{w^4s^2}\{(z+1)ws^2\frac{\partial}{\partial w} + 
 ((z + 1)^2s^4 + (z + 1)^2s^3 + w^4)\frac{\partial}{\partial s}\}\\
    & &\quad \quad   \mbox{on}~U_{\infty}\times V_{\infty}
\end{array}
$$
Since $\frac{\partial y}{\partial x}  = x^2$, we have $D^2 = x^2D$.
Therefore, the rational vector field $D$ is 2-closed.
The isolated singularities of $D$ on each affine chart are as follows.
$$
\begin{array}{ll}
  \mbox{On}~ U_0\times V_0 &   P: (x, y, t) = (0, 0, 0)\\
 \mbox{On}~ U_{\infty}\times V_{0} &   Q_1: (w, z, t) = (0, 0, 1)\\
 \mbox{On}~ U_{0} \times V_{\infty} &  \mbox{No isolated singular point}\\
 \mbox{On}~ U_{\infty} \times V_{\infty} &   
 R: (w, z, s) = (0, 0, 0), Q_2: (w, z, s) = (0, 0, 1).
 \end{array}
$$
On the surface $S$, $Q_1$ and $Q_2$ give the same point, and we denote it
by $Q$.
We set 
$$
   \psi (P) = P', \psi (Q) = \psi (R) = Q', 
   \psi^{-1}(P') = F_0, \psi^{-1}(Q') = F_{\infty}.
$$
From here on, we use the same notation for the curve and the proper transform
of the curve, if no confusion can occur. 
We blow-up at $P$, and denote the exceptional curve by $G_1$.
Then, on the exceptional curve $G_1$ there exists one isolated singular point
of the rational vector field $D$. We blow-up at the singular point, and
denote the exceptional curve by $G_2$. Then, the vector field
has no isolated singular point on $G_2$.
Now, we blow-up at $Q$, and denote the exceptional curve by $E_1$.
Then, the vector field has no isolated singular point on $E_1$.
We again blow-up at $R$, and denote the exceptional curve by $E_2$.
On the surface $\tilde{S}$ which we got by these blowing-ups
the rational vector field $D$ has no isolated singularities.
We have the morphism
$$
   \tilde{\psi} : \tilde{S} \longrightarrow E
$$
which is induced by $\psi$. Then, on $\tilde{S}$, by our construction
we have the following lemma.
\begin{lemma}\label{curves}
On $\tilde{S}$, we have the following results.

(1) $\tilde{\psi}^{-1}(P') = F_0 + G_1 + 2G_2$, 
$\tilde{\psi}^{-1}(Q') = F_{\infty} + E_1 + E_2$.

(2) The curves $F_0$, $G_1$ and $F_{\infty}$ are integral with respect to the vector field $D$. The curves
$G_2$, $E_1$, $E_2$ and $C_{\infty}$ are non-integral with respect to the vector field $D$. 

(3) $F_0^2= -2$, $G_1^2 = -2$, $G_2^2 = -1$, $F_{\infty}^2 = -2$, $E_1^2= -1$, $E_2^2 = -1$.

(4) $(F_0, G_2) = (G_2, G_1) = 1$, $(F_0, G_1) = 0$.

(5) $(F_{\infty}, E_1) = (F_1, E_2) = (C_{\infty}, E_2) = 1$, 
$(F_{\infty}, C_{\infty}) = (E_1, E_2) = (C_{\infty}, E_1) = 0$.

(6) There is a cuspidal point of the vector field $D$ on $G_2$. There is
also a cuspidal point of the vector field $D$ on $E_2$ where
it intersects with $C_{\infty}$.

\end{lemma}

We consider the quotient surface $\tilde{S}^{D}$ of $\tilde{S}$ by $D$.
We have the projection
$$
\pi : \tilde{S} \longrightarrow \tilde{S}^{D}
$$
and a commutative diagram
$$
\begin{array}{ccc}
     \tilde{S} & \stackrel{\pi}{\longrightarrow} & \tilde{S}^{D}\\
     \tilde{\psi} \downarrow &   & \downarrow \psi'  \\
       E   & \stackrel{F}{\longrightarrow} & E^{(2)}.
\end{array}
$$
Here, $F$ is the Frobenius morphism and $E^{(2)}$ is the Frobenius image.
We set $B = E^{(2)}$, $P'' = F(P')$ and $Q'' = F(Q')$.
For a curve $C$ on $\tilde{S}$, we denote the curve $\pi (C)$ on $\tilde{S}^D$
again by $C$, if no confusion can occur. 
By Lemma \ref{curves} and Proposition \ref{cuspidal curve}, 
we have the following lemma.

\begin{lemma}
On $\tilde{S}^D$, we have the following results.

(1) ${\psi'}^{-1}(P'') = 2F_0 + 2G_1 + 2G_2$, 
${\psi'}^{-1}(Q'') = 2F_{\infty} + E_1 + E_2$.

(2) $F_0^2= -1$, $G_1^2 = -1$, $G_2^2 = -2$, $F_{\infty}^2 = -1$, 
$E_1^2= -2$, $E_2^2 = -2$, $C_{\infty}^2 = -2$.

(3) $(F_0, G_2) = (G_2, G_1) = 1$, $(F_0, G_1) = 0$.

(4) $(F_{\infty}, E_1) = (F_1, E_2) = 1$, $(C_{\infty}, E_2) = 2$,
$(F_{\infty}, C_{\infty}) = (E_1, E_2) = (C_{\infty}, E_1) = 0$.

(5) $G_2$ and $E_2$ are rational cuspidal curves.
\end{lemma}

First, we blow-down $F_0$, $G_1$ and $F_{\infty}$, and 
then $E_1$ becomes an exceptional curve of the first kind and so
we blow-down it:
$$
\eta : \tilde{S}^D \longrightarrow X.
$$
Then, we have a quasi-elliptic surface
$$
  \varphi :  X \longrightarrow B.
$$
The fiber $\varphi^{-1}(P'')$ is the only one multiple fiber, and we have
no other singular fiber.

Now, let's caluculate the canonical divisor $K_X$.
First, we have
$$
    K_{\tilde{S}^D}\sim  \eta^{*}K_X + F_0 + G_1 + E_1 + 2F_{\infty}.
$$
Therefore, we have 
$$
    \pi^{*}K_{\tilde{S}^D}\sim  \pi^{*}\eta^{*}K_X + F_0 + G_1 + 2E_1 + 2F_{\infty}.
$$
On $\tilde{S}$, by a direct calculation of $D$ and $K_{\tilde{S}}$ , we have
$$
\begin{array}{l}
(D) = -2C_{\infty} -4F_{\infty} + G_1 + 4G_2 -3E_1 -3E_2,\\
K_{\tilde{S}} \sim -2C_{\infty}+ G_1 + 2G_2 + E_1 -E_2.
\end{array}
$$
Putting these data in the canonical bundle formula by Rudakov-Shafarevich:
$$
 K_{\tilde{S}} \sim (D) + \pi^{*}K_{\tilde{S}^D},
$$
we have
$$
\pi^{*}\eta^{*}K_X \sim 2(F_{\infty} + E_1 + E_2) - (F_0 + G_1 + 2G_2).
$$
Therefore, we have
$$
\eta^{*}K_X \sim (2F_{\infty} + E_1 + E_2) - (F_0 + G_1 + G_2).
$$
Hence, we have 
$$
K_X \sim E_2 - G_2  \approx G_2,
$$
where $\approx$ means numerical equivalence.
This means that there exists a divisor $\calL$ on $B$ such that
\begin{equation}\label{canonical2}
K_X \sim \varphi^{*}(\calL) + G_2.
\end{equation}
Therefore, the fiber $\varphi^{-1}(P'')$ is a tame multiple fiber.

\begin{proposition}\label{example}
The surface $\varphi: X \longrightarrow B$
which we constructed above is a quasi-elliptic surface
with only one tame
multiple fiber. It has no more singular fibers and
$\chi(\calO_{X}) = 0$ holds.
The linear system $\vert 6K_X \vert$ gives the structure of the
quasi-elliptic surface, and the linear system $\vert 5K_X \vert$ 
does not give the structure of the quasi-elliptic surface.
\end{proposition}
\proof{Take a general fiber $G$. Then, we have $G^2 = 0$ and $(K_X, G) = 0$.
Therefore, by the genus formula the virtual genus of $G$ is 1.
On the other hand, $\tilde{\psi} : \tilde{S} \longrightarrow E$ is a ruled surface.
Therefore, $G$ is not an elliptic curve. This means that 
$\varphi: X \longrightarrow B$ is a quasi-elliptic surface.
By our construction, we have Betti numbers $b_1(X) = 2$ and
$b_2(X) = 2$. Therefore, the Euler number $c_2(X) = 1 -2 + 2 - 2 + 1 = 0$.
Since $K_X^2 = 0$, we have $\chi(\calO_{X}) = 0$ by Noether's formula.
Since we have 
${\rm H}^0(X, {\calO}_X(6K_X)) \cong {\rm H}^0(B, {\calO}_B(3P''))$
and the divisor $3P''$ is very ample on $B$, 
the linear system $\vert 6K_X \vert$ gives the structure of the
quasi-elliptic surface. 
Since ${\rm H}^0(X, {\calO}_X(5K_X)) \cong {\rm H}^0(B, {\calO}_B(2P''))$
and the divisor $2P''$ is not very ample on $B$, 
the linear system $\vert 5K_X \vert$ does not give the structure of the
quasi-elliptic surface.
\hfill $\Box$}

\begin{remark}
In the above, we calculate the canonical divisor $K_X$ by
the construction of our quasi-elliptic surface.
We give here one more proof for (\ref{canonical2}).
On the quasi-elliptic surface $\varphi: X \longrightarrow B$, 
the cusp locus $C_{\infty}$ is a elliptic curve and we have
$C_{\infty}^2 = -1$ by considering the structure of blow-down.
Therefore, by the genus formula, we have $(K_X, C_{\infty}) = 1$.
On the other hand, by the canonical bundle formula 
for the quasi-elliptic surface $X$,
we have
$$
K_X \sim \varphi^{*}(\calL) + aG_2
$$
with a line bundle ${\calL}$ on $B$ and $a = 0$ or $1$. 
Since $1 = (K_X, C_{\infty}) = 2\deg~ \calL + a$,
we conclude $a = 1$ and $\deg ~\calL = 0$, which shows (\ref{canonical2}).
\end{remark}

\section{Multicanonical systems}
Let $\varphi: S \longrightarrow B$ be a quasi-elliptic surface
over an algebraically closed field $k$ of characteristic $p > 0$.
Such a surface exists only in characteristic $p = 2$ or $3$.
In this case, the multiplicity of a multiple fiber is equal to $p$
(cf. Bombieri-Mumford \cite{BM}). We denote by $pF_{i}$ $(i = 1, \ldots,
\lambda)$ the multiple fibers. Then, the canonical divisor formula
is given by
$$
     K_{S} \sim \varphi^{*}(K_{B} - {\bf f}) + \sum_{i= 1}^{\lambda}a_{i}F_{i},
$$
where ${\bf f}$ is a divisor on $B$ and 
$- \deg ~{\bf f} = \chi(\calO_{S})  + t$ with $t =$ length of
the torsion part of $R^{1}\varphi_{*}\calO_{S}$, and $0 \leq a_{i} \leq p - 1$.
For details, see Bombieri-Mumford~\cite{BM}.

We denote by $g$ the genus of the base curve $B$. Then, we have the following
theorem.

\begin{theorem}
Assume $p = 2$. Then, for any 
quasi-elliptic surface $\varphi : S \to B$ with Kodaira dimension $\kappa (S) = 1$ 
over $k$ and for any
$m \geq 6$ $\vert mK_S \vert$ gives the unique structure of quasi-elliptic surface, 
and 6 is the best possible number. 
\end{theorem}
\proof{The method of the proof is similar to the one in
Iitaka~\cite{I} and Katsura-Ueno~\cite{KU} (see also Katsura~\cite{K} \cite{K2}).
The Kodaira dimension of $S$ is equal to 1 if and only if
$$
(*) \quad 2g -2 + \chi(\calO_{S}) + t + \sum_{i=1}^{\lambda}(a_{i}/m_{i}) > 0.
$$
Therefore, we need to find the least integer $m$ such that
$$
(**)\quad m(2g -2 + \chi(\calO_{S}) + t) + \sum_{i = 1}^{\lambda}[ma_{i}/m_{i}]
\geq 2g + 1
$$
holds under the condition $(*)$. Here, $[r]$ means the integral part of 
a real number $r$. We have the following 6 cases:

Case (I) $g \geq 2$

Case (II-1) $g = 1, \chi(\calO_{S}) + t \geq 1$

Case (II-2) $g = 1, \chi(\calO_{S}) = 0, t = 0$

Case (III-1) $g = 0, \chi(\calO_{S}) + t \geq 3$

Case (III-2) $g = 0, \chi(\calO_{S}) + t = 2$

Case (III-3) $g = 0, \chi(\calO_{S})= 1, t = 0$

\noindent
Case (I) We have $2g - 2 + \chi(\calO_{S})\geq 5(g - 1)/3$.
Hence, if $m \geq 3$, $(**)$ holds.

\noindent
Case (II-1) If $m \geq 3$, $(**)$ holds.

\noindent
Case (II-2) All multiple fibers are tame in this case.
If $m \geq 6$, $(**)$ holds by $p = 2$. As we constructed in Section 4,
there exists a quasi-elliptic surface with only one 
tame multiple fiber of type II
and $\chi(\calO_{S}) = 0$ over an elliptic curve. 
Therefore, we need $m\geq 6$.

\noindent
Case (III-1) $(**)$ holds for $m \geq 1$.

\noindent
Case (III-2) Since $\chi (\calO_{S}) \geq 1$, we have $t \leq 1$. Therefore,
the number of wild fibers is less than or equal to $1$.
If there exists at least one tame multiple fiber
then $(**)$ holds for $m \geq 2$. If there exist no tame fibers and only one
wild fiber, then by Katsura-Ueno~\cite{KU} Lemma 2.4, 
this case is excluded in the case of $p = 2$.

\noindent
Case (III-3) By $p = 2$, $(**)$ holds for $m \geq 4$.

The result on the best possible number follows from the example in Section 4.
\hfill $\Box$}

\end{document}